\documentclass[11pt]{article}

\usepackage[matrix,arrow,curve,cmtip]{xy}
\usepackage{amssymb}
\usepackage{latexsym}
\usepackage{theorem}
\usepackage[a4paper,width=129mm,height=220mm]{geometry}

\def\to{\mbox{$\xymatrix@1@C=5mm{\ar@{->}[r]&}$}}
\def\tto{\mbox{$\xymatrix@1@C=5mm{\ar@{=>}[r]&}$}}
\def\halfcirc{\begin{picture}(0,0)\put(0,3){\oval(2,2)[l]}\end{picture}}
\def\surj{\mbox{$\xymatrix@1@C=5mm{\ar@{->>}[r]&}$}}
\def\inj{\mbox{$\xymatrix@1@C=5mm{\hspace{0.3mm}\ar@{>->}[r]&}$}}
\def\incl{\mbox{$\xymatrix@1@C=5mm{\ar@{->}[r]|<{\halfcirc}&}$}}
\def\bkar{\mbox{$\xymatrix@1@C=5mm{\ar@{->}[l]&}$}}
\def\distsign{\begin{picture}(0,0)\put(0,0){\circle{4}}\end{picture}}
\def\dist{\mbox{$\xymatrix@1@C=5mm{\ar@{->}[r]|{\distsign}&}$}}
\def\bkdist{\mbox{$\xymatrix@1@C=5mm{\ar@{->}[l]|{\distsign}&}$}}
\def\biar{\mbox{$\xymatrix@1@C=5mm{\ar@<1.5mm>[r]\ar@<-0.5mm>[r]&}$}}
\def\bidist{\mbox{$\xymatrix@1@C=5mm{\ar@<1.5mm>[r]|{\distsign}\ar@<-0.5mm>[r]|{\distsign}&}$}}
\def\adjar{\mbox{$\xymatrix@1@C=5mm{\ar@<1.5mm>@{<-}[r]\ar@<-0.5mm>[r]&}$}}
\def\adjdist{\mbox{$\xymatrix@1@C=5mm{\ar@<1.5mm>@{<-}[r]|{\distsign}\ar@<-0.5mm>[r]|{\distsign}&}$}}
\def\iso{\mbox{$\xymatrix@1@C=6mm{\ar@{->}[r]^{\sim}&}$}}
\def\doubiso{\mbox{$\xymatrix@1@C=6mm{\ar@{<->}[r]^{\sim}&}$}}
\def\doubar{\mbox{$\xymatrix@1@C=6mm{\ar@{<->}[r]&}$}}
\def\endoar#1#2{\mbox{\xymatrix{{#1}\ar@(u,r)|{#2}}}}
\def\ddown{\begin{picture}(6,0)
\put(0,12){\xymatrix@R=4mm{\mbox{} \\ \ar@{<<-}[u]}}
\end{picture}}
\def\down{\begin{picture}(6,0)
\put(0,12){\xymatrix@R=4mm{\mbox{} \\ \ar@{<-}[u]}}
\end{picture}}

\newtheorem{theorem}{Theorem}[section]
\newtheorem{lemma}[theorem]{Lemma}
 
\newtheorem{proposition}[theorem]{Proposition}

{\theorembodyfont{\upshape}}
\newcommand{\proof}{\noindent {\it Proof\ }: }
\def\endofproof{$\mbox{ }\hfill\Box$\par\vspace{1.8mm}\noindent}

\def\c{_{\sf c}}
\def\<{\langle}
\def\>{\rangle}
\def\sup{{\sf sup}}
\def\f{^{\sf f}}

\def\R{{\cal R}}
\def\RSDist{{\sf RSDist}}

\def\etal{~{\it et~al.}}

\def\Cocont{{\sf Cocont}}

\def\Dwn{{\sf Dwn}}
\def\Rel{{\sf Rel}}

\def\:{\colon}
\def\1{{\bf 1}}
\def\2{{\bf 2}}
\def\3{{\bf 3}}

\def\Set{{\sf Set}}
\def\QUANT{{\sf QUANT}}

\def\cc{_{\sf cc}}
\def\op{^{\sf op}}

\def\Sup{{\sf Sup}}
\def\Cat{{\sf Cat}}
\def\Matr{{\sf Matr}}
\def\Dist{{\sf Dist}}

\def\Idm{{\sf Idm}}
\def\Q{{\cal Q}}

\def\K{{\cal K}}

\def\P{{\cal P}}
\def\V{{\cal V}}

\def\U{{\cal U}}

\def\M{{\cal M}}
\def\F{{\cal F}}
\def\L{{\cal L}}

\def\colim{\mathop{\rm colim}}
\def\lim{\mathop{\rm lim}}
\def\bbA{\mathbb{A}}
\def\bbB{\mathbb{B}}
\def\bbC{\mathbb{C}}
\def\bbD{\mathbb{D}}

\def\bbP{\mathbb{P}}

\def\bbK{\mathbb{K}}

\def\bbX{\mathbb{X}}
\def\bbY{\mathbb{Y}}
\def\tensor{\otimes}

\title{Towards ``dynamic domains'': \\ totally continuous cocomplete $\Q$-categories} 
\author{Isar Stubbe\footnote{Centro de Matem\'atica, Universidade de Coimbra, Apartado 3008, 3001--454 Coimbra, Portugal, {\tt isar@mat.uc.pt}.}}
\date{October 8, 2005\footnote{This paper was presented at the 21st anual conference on Mathematical Foundations of Programming Semantics (MFPS XXI) at the University of Birmingham, 17--21 May 2005. An extended abstract was published [Stubbe, 2006].}}

\begin{document}

\maketitle

\begin{quote}
\textbf{Abstract.} It is common practice in both theoretical computer science and theoretical physics to describe the (static) logic of a system by means of a complete lattice. When formalizing the dynamics of such a system, the updates of that system organize themselves quite naturally in a quantale, or more generally, a quantaloid. In fact, we are lead to consider cocomplete quantaloid-enriched categories as fundamental mathematical structure for a dynamic logic common to both computer science and physics. Here we explain the theory of totally continuous cocomplete categories as generalization of the well-known theory of totally continuous suplattices. That is to say, we undertake some first steps towards a theory of ``dynamic domains''. 
\end{quote}

\section{Introduction}\label{A}

\paragraph{Towards ``dynamic domains''.} It is common practice in both theoretical computer science and theoretical physics to describe the `properties' of a `system' by means of a complete lattice $\L$; this lattice is then thought of as the logic of the system. For example, the lattice of closed subspaces of a Hilbert space is the logic of properties of a quantum system; and, in computer science, a domain is the logics of observables of a computational system. 
\par
More recently, also another ordered structure has been recognized to play an important r\^ole in both physics and computer science: when formalizing the dynamics of a physical or computational system, it turns out that the `updates' of a system -- think of them as programs for a computational system, and property transitions for a physical system -- organize themselves quite naturally in a quantale $\Q$ [Abramsky and Vickers, 1993; Coecke and Stubbe, 1999]. 
\par
Having a complete lattice $\L$ of properties of a system and a quantale $\Q$ of updates, we give an operational meaning to each $f\in\Q$ by the so-called Principle of Causal Duality (explained in detail in [Stubbe, 2002] but going back to [Floyd, 1967; Hoare, 1969] for computational systems and [Coecke, Moore and Stubbe, 2001] for physical systems): we want every $f\in\Q$ to determine an adjoint pair of order-preserving morphisms $f^*\dashv f_*\:\L\adjar\L$. So the left adjoint assigns to a given input $a\in\L$ its strongest consequence $f^*(a)\in\L$ under the action of $f$ (`strongest postcondition'), and the right adjoint assigns to a given output $b\in\L$ the weakest cause $f_*(b)\in\L$ under the action of $f$ (`weakest precondition'). Moreover we ask that $(g\circ f)^*=g^*\circ f^*$, $1^*=1_{\L}$ and $(\bigvee_if_i)^*=\bigvee_if_i^*$ for every $f,g,(f_i)_i\in\Q$ (and $1\in\Q$ is the unit for the monoid structure of $\Q$).
\par
In fact, a complete lattice $\L$ and a quantale $\Q$ linked by the Principle of Causal Duality, tangle up in one simple mathematical structure: a cocomplete $\Q$-enriched category. Indeed, putting $\bbA_0=\L$ as set of objects, the mapping
$$\bbA(-,-)\:\bbA_0\times\bbA_0\to\Q\:(a,b)\mapsto\bigvee\{f\in\Q\mid f^*(a)\leq b\}$$
endowes $\bbA_0=\L$ with a ``$\Q$-valued implication'' [Lawvere, 1973]: for $a,b\in\bbA_0=\L$, the element $\bbA(a,b)\in\Q$ is the weakest (i.e.\ least deterministic) update that, for input $a$, guarantees output $b$. This in fact turns $\bbA$ into a $\Q$-enriched category. This $\Q$-category is tensored and cotensored due to the Principle of Causal Duality; and the underlying order of this $\Q$-category $\bbA$ being a suplattice, namely $\L$, implies together with the tensors and cotensors that $\bbA$ is cocomplete.
\par
So, conclusively, we are lead to consider cocomplete $\Q$-categories as crucial mathematical structure in a {\em dynamic logic} as common mathematical foundation for dynamic phenomena in both computer science and physics. We will allow $\Q$ to be a quantaloid rather than a quantale, for this extra generality (allowing a `typed dynamics') doesn't really complicate matters---even though one has to bring in some adjustments to pass from enrichement in a monoidal category (i.e.\ bicategory with one object) to enrichment in a bicategory (with possibly many objects). For the basic theory of $\Q$-enriched categorical structures, see [Stubbe 2004, 2005a, 2005b]; we keep all the notations introduced there. Those works contain the more ``historical'' references on the theory of quantaloid-enriched categories.
\par
Our notation for the 2-category of $\Q$-categories and functors is $\Cat(\Q)$; and further on $\Cocont(\Q)$ denotes the 2-category of cocomplete $\Q$-categories and cocontinuous functors. 

\paragraph{Modules or cocomplete categories?}

There is an alternative and probably better known way of coupling a complete lattice $\L$ (static properties of some system) with a quantale $\Q$ (dynamics of that system): namely, by means of an action of the latter on the former. Such is a morphism $\alpha\:\L\tensor\Q\to\L$ in $\Sup$, the category of suplattices and supmorphisms (i.e.\ complete lattices and mappings that preserve arbitrary suprema), satisfying axioms on the compatibility with the monoid structure of $\Q$. Then $\L$ is said to be a (right) $\Q$-module, and with the obvious notion of homomorphism between such modules over a fixed $\Q$, one obtains a (2-)category of $\Q$-modules---which, however, is (bi)equivalent to the (2-)category $\Cocont(\Q)$ of cocomplete $\Q$-categories! 
\par
This equivalence (which also holds in the more general case of a quantaloid $\Q$, see [Stubbe, 2004] for the details) is easily sketched: regarding an action $\alpha\:\L\tensor\Q\to\L$ as a mapping $\alpha\:\L\times\Q\to\L$ that preserves suprema in both variables, it follows that for every $f\in\Q$, 
$$\alpha(-,f)\:\L\to\L\:x\mapsto\alpha(x,f)$$
has a right adjoint. So indeed, to every $f\in\Q$ we can associate an adjunction $f^*\dashv f_*$ by putting $f^*=\alpha(-,f)$. The compatibility axioms on the action then assure that the Principle of Causal Duality holds, so that -- precisely as before -- we obtain a cocomplete $\Q$-category $\bbA$ by putting $\bbA_0=\L$ as set of objects and 
$$\bbA(-,-)\:\bbA_0\times\bbA_0\to\Q\:(a,b)\mapsto\bigvee\{f\in\Q\mid\alpha(a,f)\leq b\}$$
as $\Q$-valued implication. 
\par
Conversely, given a cocomplete $\Q$-category $\bbA$, one can order its objects by the clause $a\leq b\iff 1\leq\bbA(a,b)$ (read: $a$ is smaller than $b$ when, for the system having initial property $a$, property $b$ holds after having performed the identity update) and it can be shown that this order is complete; so $\bbA_0$ is a complete lattice. Putting now $\L=\bbA_0$, the mapping
$$\L\times\Q\to\L\:(a,f)\mapsto\bigwedge\{b\in\bbA_0\mid f\leq\bbA(a,b)\}$$
which assigns to an input and an update the strongest possible output, can be shown to preserve suprema in both variables, so that it corresponds to a supmorphism $\L\tensor\Q\to\L$, which in turn proves to be an action in $\Sup$.
\par
Abramsky and Vickers [1993] (but see also [Resende, 2000] for a survey) apply the theory of $\Q$-modules to process semantics: taking into account that an informatic system may be affected by the way in which it is observed, they argue that the observable properties of an informatic system form a quantale (or even a quantaloid), and a module is then viewed as a generalization of a labelled transition system. Also in [Baltag {\it et al.}, 2004], modules on a quantale are used to cope with dynamic phenoma in computer science, in particular, to provide an algebraic semanctics for epistemic actions and updates.
\par
Our choice to work with cocomplete $\Q$-enriched categories rather than $\Q$-modules, even though they are mathematically equivalent structures, reflects a simple yet powerful idea: we explicitly put ourselves in the context of a {\em logic with truth values in $\Q$} within which we develop our mathematics. The claim in this paper is then that, even in this universe of discourse gouverned by such a ``dynamic logic'', it is possible to develop (a strong variant of) domain theory. And it is precisely because we have chosen to work with cocomplete $\Q$-categories instead of $\Q$-modules, that our presentation is so naturally a generalization of the (``classical'') results. (In section \ref{I} we shall discuss the meaning of our results for module theory though.)

\paragraph{Totally continuous suplattices.}

Suplattices are of course examples of cocomplete quantaloid-enriched categories: consider the two-element Boolean algebra $\2$ as a one-object quantaloid, then $\Sup$ is (biequivalent to) $\Cocont(\2)$. That is to say, suplattices are dynamic logics... with a trivial dynamics! Given the importance of totally continuous suplattices in computer science (as a particular kind of domain), it is natural to ask in how far the ``classical'' theory of totally continuous suplattices generalizes to $\Cocont(\Q)$. This presentation is all about giving an answer to that question. So let us first quickly recall the basics of the theory of totally continuous suplattices.
\par
On any suplattice $L$ one may define the so-called ``way-below'' relation: say that $a$ is way-below $b$, and write $a\ll b$, when for every {\em directed} downset $D\subseteq L$, $b\leq\bigvee D$ implies $a\in D$. A suplattice is said to be continuous when every element is the supremum of all elements way-below it. The theory of continuous suplattices has connections with topology and analysis (as the adjective ``continuous'' would suggest), and applications in computer science (since they are examples of ``domains''). The classical reference is [Gierz\etal, 1980].
\par
As a (stronger) variant of the above, one may also define the ``totally-below'' relation on a suplattice $L$: say that $a$ is totally-below $b$, and write $a\lll b$, when for {\em any} downset $D\subseteq L$, $b\leq\bigvee D$ implies $a\in D$. Of course $L$ is now said to be totally continuous when every element is the supremum of all elements totally-below it; in this case $L$ is also continuous. Our main reference on this subject is [Rosebrugh and Wood, 1994]. Let us recall some of the features of these structures.
\par
(a) A suplattice $L$ is totally continuous if and only if any supmorphism $f\:L\to M$ factors through any surjective supmorphism $g\:K\surj M$. This gives the totally continuous suplattices a universal status within the quantaloid $\Sup$: they are precisely its projective objects.
\par
(b) Totally continuous suplattices are precisely those suplattices for which the map sending a downset to its supremum has a left adjoint: the left adjoint to $\bigvee\:\Dwn(L)\to L\:D\mapsto\bigvee D$ is namely the map $a\mapsto\{x\in L\mid x\lll a\}$. In other words, the supremum-map is required to preserve all infima; and so such a suplattice is also said to be completely distributive\footnote{[Rosebrugh and Wood, 1994] study precisely this notion under the name of {\em constructive} complete distributivity for suplattices in a topos $\cal E$. [Fawcett and Wood, 1990] prove that, when working with suplattices in $\Set$ (and thus disposing of the axiom of choice), this constructive complete distributivity coincides with complete distributivity in the usual sense of the word. See also [Wood, 2004].}. 
\par
(c) The totally-below relation on a totally continuous suplattice is idempotent. Conversely, given a set equipped with an idempotent binary relation $(X,\prec)$, the subsets $S\subseteq X$ such that $x\in S$ if and only if there exists a $y\in S$ such that $x\prec y$, form a totally continuous suplattice. This correspondence underlies the 2-equivalence of the split-idempotent completion of $\Rel$ (whose objects are thus idempotent relations) and the full subcategory of $\Sup$ determined by the totally continuous suplattices. 
\par
(d) Given any ordered set $(X,\leq)$, the construction in (c) implies that $\Dwn(X)$ is a totally continuous suplattice. But it distinguishes itself in that every element of $\Dwn(X)$ is the supremum of ``totally compact elements'', i.e.\ elements that are totally below themselves. Such a suplattice is said to be totally algebraic; and in fact all totally algebraic suplattices are of the form $\Dwn(X)$ for some ordered set $(X,\leq)$. This correspondence underlies the 2-equivalence of the split-monad completion of $\Rel$ (whose objects are thus orders) and the full subcategory of $\Sup$ determined by the totally algebraic suplattices.
 
\paragraph{Totally continuous cocomplete $\Q$-categories.} 

In how far does the ``classical'' theory of totally continuous suplattices generalize to $\Cocont(\Q)$, the category of cocomplete $\Q$-enriched categories? The following answer is a combination of \ref{8}, \ref{103}, \ref{1000} and \ref{404} below.
\begin{theorem}
For a cocomplete $\Q$-category $\bbA$, the following are equivalent:
\begin{enumerate}
\item $\bbA$ is projective in $\Cocont(\Q)$,
\item $\bbA$ is completely distributive,
\item $\bbA$ is totally continuous,
\item $\bbA\simeq\R\bbB$ for some regular $\Q$-semicategory $\bbB$.
\end{enumerate}
And, as particular case of the above, the following are equivalent:
\begin{enumerate}
\item $\bbA$ is totally algebraic,
\item $\bbA\simeq\P\bbC$ for some $\Q$-category $\bbC$.
\end{enumerate}
Therefore, denoting $\Cocont_{\sf tc}(\Q)$, respectively $\Cocont_{\sf ta}(\Q)$, for the full sub-2-cat\-e\-go\-ry of $\Cocont(\Q)$ determined by its totally continuous objects, respectively totally algebraic objects, the following diagram, in which the horizontal equalities are bi\-equiv\-a\-lences (corestrictions of the local equivalences encountered in (\ref{1.3}) and (\ref{2}) further on), and the vertical arrows are full 2-inclusions, commutes:
$$\xymatrix@=15mm{
\RSDist(\Q)\ar@{=}[r] & \Cocont_{\sf tc}(\Q) \\
\Dist(\Q)\ar@{=}[r]\ar[u] & \Cocont_{\sf ta}(\Q)\ar[u] }$$
\end{theorem}
That is to say, the crucial aspects of the theory of totally continuous suplattices recalled above all generalize neatly to cocomplete $\Q$-categories: it is possible to make sense of such notions as `projectivity', `complete distributivity', `total continuity' and `total algebraicity' in the context of cocomplete $\Q$-categories.
\par
In the context of theoretical computer science, [Abramsky and Jung, 1994] argue that a mathematical structure deserves to be called a ``domain'' when it is an algebraic structure that unites aspects of convergence and of approximation. A totally continuous cocomplete $\Q$-category does exactly that: it is cocomplete (``every presheaf converges'') and is equipped with a well-behaved totally-below relation (``approximations from below''). The above results may then be ``translated'' into the domain theoretic lingo. For example, in section \ref{F} domain theorists will recognize the construction of bases: \ref{1000} could be read as saying that `` a cocomplete $\bbA$ is a domain if and only if it has a basis $\bbB$''. So this work really has the flavour of ``quantaloid-enriched domain theory''---or ``dynamic domains''.

\paragraph{Related work and future projects.} 

Clearly, totally continuous cocomplete $\Q$-categories are very strong structures; in particular can one argue that, having abandonned the notion of ``directedness'', their usefulness in computation is rather limited. So it is definitely an interesting project to investigate how a notion of ``directedness'' can be brought back in again. Certainly, other categorical generalizations of domain theory, in particular [Ad{\'a}mek and Rosick{\'y}, 1994; Ad{\'a}mek, 1997], may be very inspiring; our difficulty here, however, is that we need to generalize a notion such as ``directed (or filtered) colimit'' to the case of categories enriched in a quantaloid. (But it seems that Gordon and Power [1997] and also Kelly and Schmitt [2005] have ideas on that subject that will get us on track.) By the way, remark that -- precisely because we have chosen to work with the formalism of cocomplete $\Q$-categories rather than $\Q$-modules -- we have a lot of ideas and techniques from (enriched) category theory that we can try to adapt to the situation at hand!
\par
Another closely related, but at the same time very different work, is that of Wagner [1997]. Indeed, he unifies notions of ``liminf convergence'' in orders and metric spaces -- and thus gives one setting for treating recursive domain equations by a generalized inverse limit theorem \`a la Scott -- by means of categories enriched in a quantale. However, this base quantale is supposed to be {\em commutative} and its \emph{top element} is supposed to be the \emph{unit} for its multiplication. These very strong assumptions, especially the commutativity, are precisely what we want to avoid in our work: for we believe that it is an essential feature of a ``dynamic logic'' that its truth values (the possible updates of a system that constitute its dynamics) do not commute! 

\paragraph{Overview of contents.}

In section \ref{B} we first go through some considerations on monomorphisms and epimorphisms in $\Cocont(\Q)$, and show in particular that every epimorphic cocontinuous functor between cocomplete $\Q$-categories is regular. Then, in section \ref{C}, we study the `projective objects' in $\Cocont(\Q)$: we find the expected result that a projective object is precisely the retract of a $\Q$-category of presheaves. More precisely, we find that a cocomplete $\Q$-category $\bbA$ is projective if and only if the left adjoint to the Yoneda embedding $Y_{\bbA}\:\bbA\to\P\bbA$, which we denote $\sup_{\bbA}\:\P\bbA\to\bbA$ and which is an epimorphism in $\Cocont(\Q)$, admits a cocontinuous section.
\par
`Complete distributivity' is defined and studied in section \ref{D}: it is almost immediate that, for a cocomplete $\Q$-category $\bbA$, a cocontinuous section and a left adjoint to $\sup_{\bbA}\:\P\bbA\to\bbA$ are the same thing; in other words, `projectivity' and `complete distributivity' are equivalent. More involved are the results in section \ref{E}, where first, for a cocomplete $\Q$-category $\bbA$, the `totally-below relation' $\Theta_{\bbA}\:\bbA\dist\bbA$ is defined as the right extension of $\bbA(-,\sup_{\bbA}-)$ through $\P\bbA(Y_{\bbA}-,-)$ in $\Dist(\Q)$; then $\bbA$ is defined to be `totally continuous' whenever the $\Theta_{\bbA}$-weighted colimit of $1_{\bbA}$ is $1_{\bbA}$; and finally it is shown that complete distributivity and total continuity are equivalent.
\par
If $\bbA$ is totally continuous, then the totally-below relation $\Theta_{\bbA}$ is a comonad (its comultiplication is often referred to as the ``interpolation property''), and therefore an idempotent, in $\Dist(\Q)$. All idempotents split in $\RSDist(\Q)$, and the consequences thereof for the totally-below relation on a totally continuous cocomplete $\Q$-category $\bbA$ are investigated in section \ref{F}. It turns out that a cocomplete $\bbA$ is totally continuous if and only if it is (equivalent to) the category of regular presheaves on some regular $\Q$-semicategory $\bbB$.
\par
Section \ref{G} contains a discussion of so-called `totally compact objects' in a cocomplete $\Q$-category $\bbA$. Denoting $i\:\bbA\c\to\bbA$ the full embedding of $\bbA$'s totally compact objects, we define $\Sigma_{\bbA}\:\bbA\dist\bbA$ to be $\bbA(-,i-)\tensor\bbA(i-,-)$; then $\bbA$ is said to be `totally algebraic' when the $\Sigma_{\bbA}$-weighted colimit of $1_{\bbA}$ is $1_{\bbA}$. Alternatively, $\bbA$ is totally algebraic if and only if the left Kan extension of $i\:\bbA\c\to\bbA$ along itself is the identity on $\bbA$. In fact, the totally algebraic cocomplete $\Q$-categories are precisely the categories of presheaves.
\par
In section \ref{H} we briefly discuss the relation between totally algebraic cocomplete $\Q$-categories and Cauchy completions of $\Q$-categories. Finally, in section \ref{I} we consider the biequivalence of cocomplete $\Q$-categories and $\Q$-modules, and show in particular that projective modules and small-projective modules are the same thing because both these notions come down to taking retracts of direct sums of representable modules.

\paragraph{Acknowledgment.} 

The research presented in this paper was done while I was employed as Teaching and Research Assistant at the Universit\'e de Louvain (Louvain-la-Neuve, Belgium); on many occasions I had illuminating discussions with Francis Borceux, Mathieu Dupont and Jean-Roger Roisin. Further I thank Samson Abramsky, Bob Coecke and Peter Collins for inviting me to their seminars to talk about ``dynamic domains'', and for their most encouraging feedback. I am also grateful to Albert Burroni, Vincent Danos and Paul-Andre Melli\`es for the discussions we had on this matter during my visit to their \'equipe. Finally I wish to mention that my presentation of this paper at the 21st anual conference on Mathematical Foundations of Programming Semantics (MFPS XXI) at the University of Birmingham, 17--21 May 2005, wouldn't have been possible without the financial support of the organizers of that conference.

\section{Monomorphisms and epimorphisms}\label{B}

Every functor $F\:\bbA\to\bbB$ between $\Q$-categories induces an adjoint pair of distributors: $\bbB(-,F-)\:\bbA\dist\bbB$ is left adjoint to $\bbB(F-,-)\:\bbB\dist\bbA$. Now $F$ is {\em fully faithful} when the unit of this equivalence is an isomorphism, and $F$ is {\em dense} when the counit is an isomorphism. Further, the notions of a functor which is {\em essentially surjective on objects} or {\em essentially injective on objects}, speak for themselves.
\par
The locally ordered category of all (small) $\Q$-categories and functors is denoted $\Cat(\Q)$. The local order is in general not anti-symmetric so there may be non-identical isomorphic functors between two given $\Q$-categories. But an eventual isomorphism between functors is unique, and so we allow a slight abuse of language: when we say that ``the functor $F\:\bbA\to\bbB$ between $\Q$-categories is an epimorphism'', then we mean that for any $G,H\:\bbB\biar\bbC$, $G\circ F\cong H\circ F$ implies $G\cong H$; when we say that ``the functor $F\:\bbA\to\bbB$ factors through the functor $G\:\bbC\to\bbB$'', then we mean that there exists a functor $H\:\bbA\to\bbC$ such that $G\circ H\cong F$; and so on.
\par
The locally ordered category of cocomplete $\Q$-categories and cocontinuous functors is denoted $\Cocont(\Q)$. The forgetful functor $\U\:\Cocont(\Q)\to\Cat(\Q)$ admits a left adjoint (more on this in section \ref{C}), so it preserves monomorphisms. This makes the following result trivial.
\begin{lemma}\label{0.1}
For an arrow $F:\bbA\to\bbB$ in $\Cocont(\Q)$, the following are equivalent:
\begin{enumerate}
\item $F$ is a monomorphism in $\Cocont(\Q)$,
\item $F$ is a monomorphism in $\Cat(\Q)$.
\end{enumerate}
\end{lemma}
An $F\:\bbA\to\bbB$ in $\Cocont(\Q)$ has a right adjoint in $\Cat(\Q)$, say $G\:\bbB\to\bbA$. ``Taking opposites'' gives $G\op\:\bbB\op\to\bbA\op$ in $\Cocont(\Q\op)$: it is the {\em dual} of $F$, and will be denoted $F^*\:\bbB^*\to\bbA^*$. It is then quite obvious that
$$\Cocont(\Q)\to\Cocont(\Q\op)\:\Big(F\:\bbA\to\bbB\Big)\mapsto
\Big(F^*\:\bbB^*\to\bbA^*\Big)$$
is a contravariant isomorphism of 2-categories (``which is its own inverse''), so that the following is trivial.
\begin{lemma}\label{0.3}
The following are equivalent:
\begin{enumerate}
\item $F\:\bbA\to\bbB$ is an epimorphism in $\Cocont(\Q)$,
\item $F^*\:\bbB^*\to\bbA^*$ is a monomorphism in $\Cocont(\Q\op)$.
\end{enumerate}
\end{lemma}
All this now gives the following result.
\begin{proposition}\label{1}
For a left adjoint $F\:\bbA\to\bbB$ in $\Cat(\Q)$, with $F\dashv G$, the following are equivalent:
\begin{enumerate}
\item $F$ is a monomorphism in $\Cat(\Q)$,
\item $F$ is fully faithful,
\item $F$ is essentially injective on objects,
\item $G\circ F\cong 1_{\bbA}$,
\item $G$ is dense,
\item if $\bbA$ and $\bbB$ are cocomplete: $F$ is a monomorphism in $\Cocont(\Q)$.
\end{enumerate}
And also the following are equivalent:
\begin{enumerate}
\item $F$ is an epimorphism in $\Cat(\Q)$,
\item $F$ is dense,
\item $F$ is essentially surjective on objects,
\item $F\circ G\cong 1_{\bbB}$,
\item $G$ is fully faithful,
\item if $\bbA$ and $\bbB$ are cocomplete: $F$ is an epimorphism in $\Cocont(\Q)$.
\end{enumerate}
\end{proposition}
\proof
First consider the situation where $\bbA$ and $\bbB$ are not necessarily cocomplete. Since $F\circ G\circ F\cong F$, if $F$ is a monomorphism in $\Cat(\Q)$ then $G\circ F\cong 1_{\bbA}$ follows, and if $F$ is an epimorphism in $\Cat(\Q)$ then $F\circ G\cong 1_{\bbB}$ follows. All other implications are obvious and/or follow from [Stubbe, 2005a, {\bf 4.2, 4.5}].
\par
Now consider the case where both $\bbA$ and $\bbB$ are cocomplete. If $F$ is a monomorphism in $\Cocont(\Q)$, then $F$ is a monomorphism in $\Cat(\Q)$. If $F$ is an epimorphism in $\Cocont(\Q)$, then its dual $F^*=G\op$ is a monomorphism in $\Cocont(\Q\op)$, hence fully faithful in $\Cocont(\Q\op)$, so $G$ is fully faithful in $\Cocont(\Q)$. The remaining implications are trivial.
\endofproof
Part of the above is ``abstract nonsense'', i.e.\ valid in any locally ordered category and not just $\Cat(\Q)$. 
\begin{proposition}\label{0.51}
Every epimorphism in $\Cocont(\Q)$ is regular.
\end{proposition}
\proof
Let $F\:\bbA\to\bbB$ be any morphism in $\Cocont(\Q)$. It is easy to see that 
$$\bbK_0=\{(a_1,a_2)\in\bbA_0\times\bbA_0\mid Fa_1\cong Fa_2\},$$  $$\bbK\Big((b_1,b_2),(a_1,a_2)\Big)=\bbA(b_1,a_1)\wedge\bbA(b_2,a_2)$$
defines a $\Q$-category $\bbK$, and that 
$$D_1\:\bbK\to\bbA\:(a_1,a_2)\mapsto a_1,\ D_2\:\bbK\to\bbA\:(a_1,a_2)\mapsto a_2$$
are functors satisfying $F\circ D_1\cong F\circ D_2$. Now consider a weighted colimit diagram in $\bbK$, like so:
$$\xymatrix@=15mm{
\bbX\ar[r]|{\distsign}^{\Phi} & \bbC\ar[r]^{H} & \bbK.}$$
Both $\colim(\Phi,D_1\circ H)$ and $\colim(\Phi,D_2\circ H)$ exist (because $\bbA$ is cocomplete), and their images by $F$ are isomorphic (because $F$ is cocontinuous and equalizes $D_1$ and $D_2$); so certainly is $(\colim(\Phi,D_1\circ H),\colim(\Phi,D_2\circ H))$ an object of $\bbK$. But for any $(x_1,x_2)\in\bbK$ we can calculate that:
\begin{eqnarray*}
 & & \bbK\Big((\colim(\Phi,D_1\circ H),\colim(\Phi,D_2\circ H)),(x_1,x_2)\Big) \\
 & & = \bbA\Big(\colim(\Phi,D_1\circ H),x_1\Big)\wedge\bbA\Big(\colim(\Phi,D_2\circ H),x_2\Big) \\
 & & = \Big[\Phi,\bbA(D_1\circ H-,x_1)\Big]\wedge\Big[\Phi,\bbA(D_2\circ H-,x_2)\Big] \\
 & & = \Big[\Phi,\bbA(D_1\circ H-,x_1)\wedge\bbA(D_2\circ H-,x_2)\Big] \\
 & & = \Big[\Phi,\bbK(H-,(x_1,x_2))\Big].
\end{eqnarray*}
(We used that $\Phi\tensor-\dashv[\Phi,-]$ in $\Dist(\Q)$ to pass from the third to the fourth line.) That is to say, $(\colim(\Phi,D_1\circ H),\colim(\Phi,D_2\circ H))$ is the $\Phi$-weighted colimit of $H$ in $\bbK$. From this ``componentwise'' construction of colimits in $\bbK$ it immediately follows that $D_1$ and $D_2$ are cocontinuous. So we have
$$\xymatrix@=15mm{
\bbK\ar@<1mm>[r]^{D_1}\ar@<-1mm>[r]_{D_2} & \bbA\ar[r]^F & \bbB}$$
in $\Cocont(\Q)$, and $F\circ D_1\cong F\circ D_2$.
\par
If now $F\:\bbA\to\bbB$ is an epimorphism, i.e.\ $F\circ G\cong 1_{\bbB}$ for $F\dashv G$ in $\Cat(\Q)$, then we can prove that the diagram above is a universal coequalizer diagram: we claim that for a morphism $F'\:\bbA\to\bbB'$ in $\Cocont(\Q)$ such that $F'\circ D_1\cong F'\circ D_2$, the functor $\overline{F'}=F'\circ G$ is the essentially unique cocontinuous factorization of $F'$ through $F$:
$$\xymatrix@=15mm{
\bbK\ar@<1mm>[r]^{D_1}\ar@<-1mm>[r]_{D_2} & \bbA\ar@<2.5mm>@{}[r]|{\top}\ar@{->>}[r]_F\ar[rd]_{F'} & \bbB\ar@{.>}[d]^{\overline{F'}}\ar@/_5mm/[l]_{G} \\
 & & \bbB'}$$
First note that, by the assumption $F'\circ D_1\cong F'\circ D_2$ and by the construction of $\bbK$, if $Fa\cong Fa'$ then also $F'a\cong F'a'$. Now consider a weighted colimit diagram
$$\xymatrix@=15mm{
\bbY\ar[r]|{\distsign}^{\Psi} & \bbD\ar[r]^{K} & \bbB.}$$
Using the surjectivity of $F$ we have
\begin{eqnarray*}
F\circ G\circ\colim(\Psi,K)
 & \cong & \colim(\Psi,K) \\
 & \cong & \colim(\Psi,F\circ G\circ K) \\
 & \cong & F\circ\colim(\Psi,G\circ K)
\end{eqnarray*}
and therefore also
$F'\circ G\circ\colim(\Psi,K)\cong F'\circ\colim(\Psi,G\circ K)$,
from which in turn, using now the cocontinuity of $F'$,
$\overline{F'}\circ\colim(\Psi,K)\cong\colim(\Psi,\overline{F'}\circ K)$.
This shows that $\overline{F'}$ is cocontinuous. But since $F\circ G\circ F\cong F$ it follows that $F'\circ G\circ F\cong F'$, so $\overline{F'}$ is a factorization of $F'$ through $F$. This factorization is essentially unique because $F$ is an epimorphism.
\endofproof
The proof above is really a generalization of the typical direct proof of the fact that all epimorphisms in $\Sup$ are regular: recall that $\Sup\simeq\Cocont(\2)$, so we generalized the ``classical'' (i.e.\ $\2$-enriched) case to the $\Q$-enriched case. (And actually, the $D_1,D_2\:\bbK\biar\bbA$ as constructed in the first part of the proof for {\em any} $F\:\bbA\to\bbB$ in $\Cocont(\Q)$, are its kernel pair.)
\par
In what follows we will often speak of {\em surjections} in $\Cocont(\Q)$ when we mean epimorphisms.

\section{Projective cocomplete $\Q$-categories}\label{C}

The forgetful 2-functor $\U\:\Cocont(\Q)\to\Cat(\Q)$ admits a left 2-adjoint: the free cocompletion of a $\Q$-category $\bbA$ is the presheaf category $\P\bbA$. By a {\em free object} in $\Cocont(\Q)$ we will mean a free object relative to the forgetful functor $\U$, i.e.\ an object equivalent to the presheaf category $\P\bbA$ on some $\Q$-category $\bbA$.
\par
In fact, the free 2-functor $\P\:\Cat(\Q)\to\Cocont(\Q)$ is the composition of two 2-functors. First every functor $F\:\bbA\to\bbB$ induces a left adjoint distributor (the ``graph'' of $F$),
\begin{equation}\label{1.31}
\Cat(\Q)\to\Dist(\Q)\:\Big(F\:\bbA\to\bbB\Big)\mapsto\Big(\bbB(-,F-)\:\bbA\to\bbB\Big).
\end{equation}
Then every distributor determines a cocontinuous functor between presheaf categories,
\begin{equation}\label{1.3}
\Dist(\Q)\to\Cocont(\Q)\:\Big(\Phi\:\bbA\dist\bbB\Big)\mapsto\Big(\Phi\tensor -\:\P\bbA\to\P\bbB\Big).
\end{equation}
The latter is locally an equivalence (actually, locally an isomorphism since $\Dist(\Q)$ is a quantaloid and each $\P\bbB$ is skeletal). There are more details in [Stubbe, 2005a, {\bf 3.7, 6.12}].
\par
The adjunction $\P\dashv\U$ works as follows: a functor $F\:\bbA\to\bbB$ from any $\Q$-cat\-e\-go\-ry into a cocomplete $\Q$-category determines a cocontinuous functor $\<F,Y_{\bbA}\>\:\P\bbA\to\bbB$ by (pointwise) left Kan extension of $F$ along the Yoneda embedding for $\bbA$; and a cocontinuous functor $G\:\P\bbA\to\bbB$ into a cocomplete $\Q$-category determines a functor $G\circ Y_{\bbA}\:\bbA\to\bbB$ by composition with the Yoneda embedding. In other words, for an $\bbA\in\Cat(\Q)$, the Yoneda embedding $Y_{\bbA}\:\bbA\to\P\bbA$ gives the unit of the adjunction; and for some $\bbB\in\Cocont(\Q)$, the left Kan extension $\<1_{\bbB},Y_{\bbB}\>\:\P\bbB\to\bbB$ gives the counit. The latter sends a presheaf $\phi\in\P\bbB$ to the colimit $\colim(\phi,1_{\bbB})$, and will be denoted from now on as $\sup_{\bbB}\:\P\bbB\to\bbB$ (for ``supremum'' of course). Actually, $\sup_{\bbB}$ is left adjoint to $Y_{\bbB}$ in $\Cat(\Q)$; since the latter is fully faithful, the former is surjective. We refer to [Stubbe, 2005a, sections 5 and 6] for details.
\par
A {\em projective object} $\bbA$ in $\Cocont(\Q)$ is one such that in $\Cocont(\Q)$ any arrow $F\:\bbA\to\bbB$ factors (up to local isomorphism) through any surjection $G\:\bbC\surj\bbB$. This definition is classical for ordinary categories\footnote{Usually one defines ``projectivity'' with respect to a preferred class of epimorphisms, giving rise to ``regular projectivity'', ``strong projectivity'', and whatnot. But every epimorphism in $\Cocont(\Q)$ is regular, so we speak of ``projectivity'' {\em tout court}. See also section \ref{I}.}, and the following lemmas will surely ring a bell [Borceux, 1994].
\begin{lemma}\label{5}
The retract of a projective object in $\Cocont(\Q)$ is again projective.
\end{lemma}
\proof
Suppose that $S\:\bbA\to\bbP$ and $P\:\bbP\to\bbA$ exhibit $\bbA$ as retract of a projective object $\bbP$ in $\Cocont(\Q)$. Given an arrow $F\:\bbA\to\bbC$ and a surjection $G\:\bbB\surj\bbC$ in $\Cocont(\Q)$, the projectivity of $\bbP$ implies the existence of an arrow $H\:\bbP\to\bbB$ satisfying $G\circ H\cong F\circ P$, so that $H\circ S$ is a factorization of $F$ through $G$. 
\endofproof
\begin{lemma}\label{6}
Free objects in $\Cocont(\Q)$ are projective.
\end{lemma}
\proof
Consider $F\:\P\bbA\to\bbB$ and $G\:\bbC\surj\bbB$ in $\Cocont(\Q)$; so the latter is surjective. The cocontinuous functor $F\:\P\bbA\to\bbB$ corresponds, under the adjunction $\P\dashv\U$, to the functor $F\circ Y_{\bbA}\:\bbA\to\bbB$ in $\Cat(\Q)$. Denoting $H\:\bbB\to\bbC$ for the right adjoint section to $G$ in $\Cat(\Q)$ (see \ref{1}), surely the functor $H\circ F\circ Y_{\bbA}$ is a factorization of $F\circ Y_{\bbA}$ through $G$ in $\Cat(\Q)$. Again under the ``free cocompletion'' adjunction, the functor $H\circ F\circ Y_{\bbA}\:\bbA\to\bbB$ corresponds to the cocontinuous functor $\<H\circ F\circ Y_{\bbA},Y_{\bbA}\>\:\P\bbA\to\bbB$. This latter functor is a pointwise left Kan extension, hence by cocontinuity and surjectivity of $G$,
$$G\circ \<H\circ F\circ Y_{\bbA},Y_{\bbA}\>\cong \<G\circ H\circ F\circ Y_{\bbA},Y_{\bbA}\>\cong \<F\circ Y_{\bbA},Y_{\bbA}\>\cong F.$$ That is, $\<H\circ F\circ Y_{\bbA},Y_{\bbA}\>$ is a factorization of $F$ through $G$ in $\Cocont(\Q)$.
\endofproof
It follows that $\Cocont(\Q)$ {\em has enough projectives}, i.e.\ that every object in $\Cocont(\Q)$ is the quotient of a projective object: there is always the surjection $\sup_{\bbA}\:\P\bbA\surj\bbA$.
\begin{proposition}\label{7}
For a cocomplete $\Q$-category $\bbA$, the following are equivalent:
\begin{enumerate}
\item $\bbA$ is a projective object in $\Cocont(\Q)$,
\item $\sup_{\bbA}\:\P\bbA\surj\bbA$ has a section in $\Cocont(\Q)$,
\item $\bbA$ is a retract of $\P\bbA$ in $\Cocont(\Q)$,
\item $\bbA$ is a retract of a free object in $\Cocont(\Q)$.
\end{enumerate}
\end{proposition}
\proof
If $\bbA$ is a projective object in $\Cocont(\Q)$, then there must be a factorization of $1_{\bbA}\:\bbA\to\bbA$ through the surjection $\sup_{\bbA}\:\P\bbA\surj\bbA$. This proves that $\bbA$ is a retract of the free object $\P\bbA$. The remainder of the proof follows from \ref{5} and \ref{6}.
\endofproof
\par
The definition of `projective object' (in $\Cocont(\Q)$, or in any category for that matter) guarantees the existence of certain factorizations, but does not explain a way of calculating them. But in $\Cocont(\Q)$ ``liftings provide factorizations'' as soon as the latter are known to exist. (Note that right liftings always exist in $\Cocont(\Q)$: for it is a locally ordered, locally small category with stable local colimits---see [Stubbe 2005a, {\bf 6.12}].)
\begin{lemma}\label{7.1}
For $F\:\bbA\to\bbB$ and $G\:\bbC\to\bbB$ in $\Cocont(\Q)$, if a factorization of $F$ through $G$ exists, then also the right lifting $[G,F]\:\bbA\to\bbC$ of $F$ through $G$ is such a factorization.
\end{lemma}
\proof
Suppose that $H\:\bbA\to\bbB$ in $\Cocont(\Q)$ satisfies $G\circ H\cong F$. Then, by the universal property of the right lifting, $H\leq[G,F]$. But this in turn implies that $F\cong G\circ H\leq G\circ[G,F]\leq F$, so $G\circ[G,F]\cong F$.
\endofproof
For example, the factorization calculated in the proof of \ref{6} is the right lifting.
\begin{proposition}\label{7.2}
For a cocomplete $\Q$-category $\bbA$, the following are equivalent:
\begin{enumerate}
\item $\bbA$ is a projective object in $\Cocont(\Q)$,
\item for any $F\:\bbA\to\bbB$ and surjective $G\:\bbC\surj\bbB$ in $\Cocont(\Q)$, $G\circ [G,F]\cong F$.
\end{enumerate}
\end{proposition}
\par
Clearly there is a more abstract setting for these results: if, in a locally ordered category $\K$, $f\:A\to B$ factors through $g\:C\to B$ and moreover the right lifting of $f$ through $g$ exists, then the lifting is also a factorization. This presumably lead [Rosebrugh and Wood, 1994] to say that an object $A\in\K$ is {\em universally projective} when, for every $f\:A\to B$ and ``surjective'' $g\:C\surj B$, the lifting of $f$ through $g$ exists and is a factorization, i.e.\ that $g\circ [g,f]\cong f$ in $\K$. Here ``surjectivity'' must be given a meaning in $\K$; thereto [Rosebrugh and Wood, 1994] consider a proarrow equipment $(-)_{\#}\:\K\to\M$, and call $g\:C\to B$ in $\K$ ``surjective'' when the counit for the left adjoint $g_{\#}\:C\to B$ in $\M$ is an isomorphism. In those terms then, projective objects and universally projective objects are the same thing in $\Cocont(\Q)$, when considering the ``forgetful'' proarrow equipment $\Cocont(\Q)\to\Cat(\Q)$.

\section{Completely distributive cocomplete $\Q$-categories}\label{D}

A {\em (constructively\footnote{We will not insist on the adjective ``constructive'' as do [Rosebrugh and Wood, 1994], because we think that, in the context of $\Q$-categories, no confusion will arise. However, note that all our proofs are indeed constructive, which is consistent with the idea that $\Q$ is the object of truth values of a ``dynamic logic'' within which we work.}) completely distributive} cocomplete $\Q$-category $\bbA$ is one for which the left adjoint to the Yoneda embedding, $\sup_{\bbA}\:\P\bbA\surj\bbA$, has a further left adjoint. The terminology is classical for $\Q=\2$, i.e.\ for suplattices [Rosebrugh and Wood, 1994].
\begin{proposition}\label{8}
For a cocomplete $\Q$-category $\bbA$, the following are equivalent:
\begin{enumerate}
\item $\bbA$ is completely distributive,
\item $\bbA$ is a projective object in $\Cocont(\Q)$.
\end{enumerate}
\end{proposition}
\proof
Suppose that $L\dashv\sup_{\bbA}$ in $\Cat(\Q)$. Then $L$ is cocontinuous (because it is a left adjoint) and fully faithful (because $\sup_{\bbA}$ is surjective), so $\sup_{\bbA}\circ L\cong 1_{\bbA}$. That is to say, $L$ is a section to $\sup_{\bbA}$ in $\Cocont(\Q)$. Conversely, if $S\:\bbA\to\P\bbA$ is a cocontinuous section to $\sup_{\bbA}\:\P\bbA\surj\bbA$, then $\sup_{\bbA}\circ S\cong 1_{\bbA}$ implies $S\leq Y_{\bbA}$ (because $\sup_{\bbA}\dashv Y_{\bbA}$), and hence, for any $\phi\in\P\bbA$,
$$S\circ\sup_{\bbA}(\phi)\cong\colim(\phi,S)\leq\colim(\phi,Y_{\bbA})\cong\phi$$
(because $S$ is cocontinuous). So $S\circ\sup_{\bbA}\leq 1_{\P\bbA}$, which proves it to be left adjoint to $\sup_{\bbA}$.
\endofproof
The above says that, for a cocomplete $\Q$-category $\bbA$, a cocontinuous section to $\sup_{\bbA}\:\P\bbA\surj\bbA$ is the same thing as a left adjoint. But there may be several non-cocontinuous sections for $\sup_{\bbA}$, e.g.~the Yoneda embedding!
\par
Since $\Cocont(\Q)$ is a locally ordered category in which both right extensions and right liftings exist, we can use these to ``approximate'' left adjoints and cocontinuous sections to $\sup_{\bbA}\:\bbA\to\P\bbA$. Our notations are
$$\xymatrix@=15mm{
\P\bbA\ar@{->>}[d]_{\sup_{\bbA}}\ar[r]^{1_{\P\bbA}} & \P\bbA \\
\bbA\ar@{.>}[ur]_{\{\sup_{\bbA},1_{\P\bbA}\}}}\hspace{10mm}
\xymatrix@=15mm{
 & \P\bbA\ar@{->>}[d]^{\sup_{\bbA}} \\
\bbA\ar[r]_{1_{\bbA}}\ar@{.>}[ur]^{[\sup_{\bbA},1_{\bbA}]} & \bbA}$$
for the right extension of $1_{\P\bbA}$ through $\sup_{\bbA}$, respectively the right lifting of $1_{\bbA}$ through $\sup_{\bbA}$.
\begin{proposition}\label{101}
For a cocomplete $\Q$-category $\bbA$, the following are equivalent:
\begin{enumerate}
\item $\bbA$ is completely distributive,
\item $\sup_{\bbA}\circ \{\sup_{\bbA},1_{\P\bbA}\}\geq 1_{\bbA}$,
\item $\sup_{\bbA}\circ \{\sup_{\bbA},1_{\P\bbA}\}\cong 1_{\bbA}$.
\end{enumerate}
In this case, $\{\sup_{\bbA},1_{\P\bbA}\}$ is the left adjoint to $\sup_{\bbA}$ (and therefore also its cocontinuous section).
\end{proposition}
\proof
In any locally ordered category $\K$, an arrow $f\:A\to B$ has a left adjoint if and only if the right extension $\{f,1_B\}\:B\to A$ of $1_B$ through $f$ exists and satisfies $f\circ \{f,1_{B}\}\geq 1_B$; in this case, $\{f,1_B\}\dashv f$ in $\K$. Applied to $\Cocont(\Q)$, this proves the equivalence of the first and the second statement. The second and the third are equivalent because the left adjoint to $\sup_{\bbA}$ is automatically its cocontinuous section, and {\em vice versa}.
\endofproof
\begin{proposition}\label{12}
For a cocomplete $\Q$-category $\bbA$, the following are equivalent:
\begin{enumerate}
\item $\bbA$ is completely distributive,
\item $\sup_{\bbA}\circ [\sup_{\bbA},1_{\bbA}]\geq 1_{\bbA}$,
\item $\sup_{\bbA}\circ [\sup_{\bbA},1_{\bbA}]\cong 1_{\bbA}$.
\end{enumerate}
In this case, $[\sup_{\bbA},1_{\bbA}]$ is the cocontinuous section to $\sup_{\bbA}$ (and therefore also its left adjoint).
\end{proposition}
\proof
The second and third statement are always equivalent, because $\sup_{\bbA}\circ [\sup_{\bbA},1_{\bbA}]\leq 1_{\bbA}$. Since for a cocomplete $\Q$-category $\bbA$ there is at most one cocontinuous section to $\sup_{\bbA}\:\P\bbA\surj\bbA$, \ref{7.2} implies that $[\sup_{\bbA},1_{\bbA}]$ is the only candidate for the job. So if $\bbA$ is completely distributive, then $[\sup_{\bbA},1_{\bbA}]$ is the cocontinuous section to $\sup_{\bbA}$. Conversely, if $\sup_{\bbA}\circ [\sup_{\bbA},1_{\bbA}]\cong 1_{\bbA}$ then $\bbA$ is a retract of a free object, so (by \ref{8}) it is completely distributive.
\endofproof
The results in \ref{101} and \ref{12} will not be used further on.

\section{Totally continuous cocomplete $\Q$-categories}\label{E}

Given a completely distributive cocomplete $\Q$-category $\bbA$, the left adjoint to the surjection $\sup_{\bbA}\:\P\bbA\surj\bbA$ is a functor, say $T_{\bbA}\:\bbA\to\P\bbA$, satisfying 
$$\P\bbA(T_{\bbA}-,-)=\bbA(-,\sup_{\bbA}-).$$ 
By the universal property of the presheaf category $\P\bbA$, this functor -- like any functor from $\bbA$ to $\P\bbA$, for that matter -- determines, and is determined by, a distributor $\Theta_{\bbA}\:\bbA\dist\bbA$ through the formula $T_{\bbA}(a)(a')=\Theta_{\bbA}(a',a)$ [Stubbe, 2005a, {\bf 6.1}]. The elements of this distributor can be written as
\begin{eqnarray*}
\Theta_{\bbA}(a',a)
 & = & \P\bbA(Y_{\bbA}a',T_{\bbA}a) \\
 & = & \P\bbA(Y_{\bbA}a',-)\tensor\P\bbA(-,T_{\bbA}a) \\
 & = & \{\bbA(T_{\bbA}a,-),\P\bbA(Y_{\bbA}a',-)\} \\
 & = & \{\bbA(a,\sup_{\bbA}-),\P\bbA(Y_{\bbA}a',-)\}.
\end{eqnarray*}
That is to say, for a completely distributive cocomplete $\Q$-category $\bbA$ the distributor $\Theta_{\bbA}$ is the right extension of $\bbA(-,\sup_{\bbA}-)$ through $\P\bbA(Y_{\bbA}-,-)$ in $\Dist(\Q)$:
$$\xymatrix@=20mm{
\P\bbA\ar[r]|{\distsign}^{\P\bbA(Y_{\bbA}-,-)}\ar[d]|{\distsign}_{\bbA(-,\sup_{\bbA}-)}
 & \bbA \\
\bbA\ar@{.>}[ur]|{\distsign}_{\Theta_{\bbA}=\{\bbA(-,\sup_{\bbA}-),\P\bbA(Y_{\bbA}-,-)\}.}}$$
But this right extension makes sense for {\em any} cocomplete $\Q$-category $\bbA$, so -- whether $\bbA$ is completely distributive or not -- we can {\em define} the distributor $\Theta_{\bbA}\:\bbA\dist\bbA$ to be this right extension, and denote $T_{\bbA}\:\bbA\to\P\bbA$ for the functor corresponding with $\Theta_{\bbA}$ under the universal property of $\P\bbA$. In analogy with the case $\Q=\2$, we call the distributor $\Theta_{\bbA}\:\bbA\dist\bbA$ the {\em totally-below relation} on the cocomplete $\Q$-category $\bbA$; and the functor $T_{\bbA}\:\bbA\to\P\bbA$ sends an object $a\in\bbA$ to the ``presheaf of objects totally-below $a$''. The calculation rules for weighted colimits [Stubbe, 2005a, {\bf 5.2}]  make the following trivial.
\begin{lemma}\label{104}
For a cocomplete $\Q$-category $\bbA$, the following are equivalent:
\begin{enumerate}
\item for every $a\in\bbA$, $\sup_{\bbA}(T_{\bbA}a)\cong a$,
\item $\sup_{\bbA}\circ T_{\bbA}\cong 1_{\bbA}$,
\item $\colim(\Theta_{\bbA},1_{\bbA})\cong 1_{\bbA}$.
\end{enumerate}
\end{lemma}
A cocomplete $\Q$-category $\bbA$ is said to be {\em totally continuous} when it satisfies the equivalent conditions above; that is to say, ``every object in $\bbA$ is the supremum of the objects totally-below it''. We will see in \ref{103} that ``totally continuous'' is synonymous with ``completely distributive''. But first we record two helpful lemmas, the first of which literally is the ``classical'' definition of `totally-below' (when we put $\Q=\2$)!
\begin{lemma}\label{104.1}
For a cocomplete $\Q$-category $\bbA$, the elements of the totally-below relation $\Theta_{\bbA}\:\bbA\dist\bbA$ are, for $a,a'\in\bbA$,
$$\Theta_{\bbA}(a',a)=\bigwedge_{\phi\in\P\bbA}\{\bbA(a,\sup_{\bbA}\phi),\phi(a')\}.$$
\end{lemma}
\proof
By definition, $\Theta_{\bbA}$ is a right extension in $\Dist(\Q)$; with the Yoneda lemma for $\Q$-categories, an explicit calculation of this extension gives
\begin{eqnarray*}
\Theta_{\bbA}(a',a)
 & = & \{\bbA(a,\sup_{\bbA}-),\P\bbA(Y_{\bbA}a',-)\} \\
 & = & \bigwedge_{\phi\in\P\bbA}\{\bbA(a,\sup_{\bbA}\phi),\P\bbA(Y_{\bbA}a',\phi)\} \\
 & = & \bigwedge_{\phi\in\P\bbA}\{\bbA(a,\sup_{\bbA}\phi),\phi(a')\}
\end{eqnarray*}
which is precisely the claimed formula.
\endofproof
\begin{lemma}\label{103.0}
For a cocomplete $\Q$-category $\bbA$ we have that the totally-below relation $\Theta_{\bbA}\:\bbA\dist\bbA$ satisfies $\Theta_{\bbA}\leq\bbA$ and $\Theta_{\bbA}\tensor\Theta_{\bbA}\leq\Theta_{\bbA}$.
\end{lemma}
\proof
For any $a,a'\in\bbA$, put $\phi=Y_{\bbA}a$ in \ref{104.1} and use that $\sup_{\bbA}\circ Y_{\bbA}\cong 1_{\bbA}$ to calculate that $\Theta_{\bbA}(a',a)\leq\{\bbA(a,a),\bbA(a',a)\}$. Hence -- since $1_{ta}\leq\bbA(a,a)$ -- $\Theta_{\bbA}(a',a)\leq\bbA(a',a)$. The second inequality follows trivially.
\endofproof
\begin{proposition}\label{103}
For a cocomplete $\Q$-category $\bbA$, the following are equivalent:
\begin{enumerate}
\item $\bbA$ is completely distributive,
\item $\bbA$ is totally continuous.
\end{enumerate}
In this case, $T_{\bbA}$ is the left adjoint to $\sup_{\bbA}$ (and therefore also its cocontinuous section).
\end{proposition}
\proof
By \ref{103.0} the functor $T_{\bbA}\:\bbA\to\P\bbA$ satisfies $T_{\bbA}\circ\sup_{\bbA}\leq Y_{\bbA}\circ\sup_{\bbA}\leq 1_{\P\bbA}$ (whether $\bbA$ is completely distributive or not). So the second statement implies that $T_{\bbA}\dashv\sup_{\bbA}$, that is, $\bbA$ is completely distributive. Conversely, if $\bbA$ is completely distributive then, as argued in the beginning of this section, $T_{\bbA}\dashv\sup_{\bbA}$, so -- by surjectivity of $\sup_{\bbA}$ -- $\sup_{\bbA}\circ T_{\bbA}\cong 1_{\bbA}$. 
\endofproof
\par
The single most important property of the totally-below relation on a (totally continuous) cocomplete $\Q$-category is the following.
\begin{proposition}\label{105}
Given a totally continuous cocomplete $\Q$-category $\bbA$, the total\-ly-below relation $\Theta_{\bbA}\:\bbA\dist\bbA$ is a comonad in $\Dist(\Q)$.
\end{proposition}
\proof
For $a\in\bbA$, consider the presheaf $\Theta_{\bbA}\tensor\Theta_{\bbA}(-,a)$ on $\bbA$; by the calculation rules for weighted colimits [Stubbe, 2005a, {\bf 5.2}] and the result in \ref{104},
\begin{eqnarray*}
\sup_{\bbA}\Big(\Theta_{\bbA}\tensor\Theta_{\bbA}(-,a)\Big)
 & \cong & \colim\Big(\Theta_{\bbA}(-,a),\colim(\Theta_{\bbA},1_{\bbA})\Big) \\
 & \cong & \sup_{\bbA}\Big(\Theta_{\bbA}(-,a)\Big) \\
 & \cong  & a.
\end{eqnarray*}
Putting $\phi=\Theta_{\bbA}\tensor\Theta_{\bbA}(-,a)$ in \ref{104.1} gives
$$\Theta_{\bbA}(a',a)\leq\Big\{\bbA(a,a),\Theta_{\bbA}(a',-)\tensor\Theta_{\bbA}(-,a)\Big\}$$
which -- since $1_{ta}\leq\bbA(a,a)$ -- implies that $\Theta_{\bbA}(a',a)\leq\Theta_{\bbA}(a',-)\tensor\Theta_{\bbA}(-,a)$. This proves that $\Theta_{\bbA}\leq\Theta_{\bbA}\tensor\Theta_{\bbA}$, which together with \ref{103.0} gives the result.
\endofproof
The comultiplication of $\Theta_{\bbA}$ is often called its {\em interpolation property}. The result implies in particular that the totally-below relation on a totally continuous cocomplete $\Q$-category is idempotent.

\section{Splitting the totally-below relation}\label{F}

Recall from [Stubbe, 2005b, {\bf 4.5}] that, considering regular $\Q$-semicategories and regular semidistributors,
\begin{equation}\label{2}
\RSDist(\Q)\to\Cocont(\Q)\:\Big(\Phi\:\bbA\dist\bbB\Big)\mapsto\Big(\Phi\tensor -\:\R\bbA\to\R\bbB\Big)
\end{equation}
is locally an equivalence. In particular, a cocontinuous functor $F\:\R\bbA\to\R\bbB$ determines (and is determined by) the regular semidistributor $\Phi\:\bbA\dist\bbB$ with elements $\Phi(b,a)=F(Y_{\bbA}(a))(b)$. Note that $\Dist(\Q)$ is a full subquantaloid of $\RSDist(\Q)$, and that the domain restriction of (\ref{2}) to $\Dist(\Q)$ is the local equivalence in (\ref{1.3}): for a $\Q$-category $\bbA$, $\R\bbA=\P\bbA$. 
\par
Furthermore, [Stubbe, 2005b, {\bf 3.12}] says that, for each regular $\Q$-semicat\-e\-gory $\bbB$, the $\Q$-category $\R\bbB$ of regular presheaves on $\bbB$ is an essential (co)locali\-za\-tion of a certain presheaf category. So certainly is $\R\bbB$ a projective object in $\Cocont(\Q)$, i.e.\ a totally continuous cocomplete $\Q$-category  (see \ref{8} and \ref{103}). In fact, all totally continuous cocomplete $\Q$-categories are of the form $\R\bbB$, for some regular $\Q$-semicategory $\bbB$, as we show next.
\begin{proposition}\label{1000}
For a cocomplete $\Q$-category $\bbA$, the following are equivalent:
\begin{enumerate}
\item $\bbA$ is totally continuous,
\item $\bbA\simeq\R\bbB$ in $\Cocont(\Q)$ for some regular $\Q$-semicategory $\bbB$.
\end{enumerate}
In this case, the ``$\bbB$'' in the second statement is the regular $\Q$-semicategory, unique up to Morita equivalence\footnote{See [Stubbe, 2005b, section 4] for a discussion of ``Morita equivalence'' for regular $\Q$-semicategories.}, over which the totally-below relation on $\bbA$, $\Theta_{\bbA}\:\bbA\dist\bbA$, splits in $\RSDist(\Q)$.
\end{proposition}
\proof
Suppose that $\bbA$ is a totally continuous cocomplete $\Q$-category. The totally-below relation $\Theta_{\bbA}\:\bbA\dist\bbA$ is an idempotent in $\Dist(\Q)$ (see \ref{105}), hence an idempotent in $\RSDist(\Q)$. But in the latter quantaloid idempotents split [Stubbe, 2005b, Appendix] so there must exist a regular $\Q$-semicategory, unique up to Morita equivalence, over which $\Theta_{\bbA}$ splits; let us denote such a splitting as
$$\xymatrix@=15mm{
\bbA\ar@(u,l)_{\Theta_{\bbA}}|{\distsign}\ar@/^3mm/[r]^{\Phi}|{\distsign} & \bbB\ar@/^3mm/[l]^{\Psi}|{\distsign}}.$$
Note that $\Psi\dashv\Phi$ (because $\Theta_{\bbA}\leq\bbA$), so that applying (\ref{2}) we may now consider the diagram
$$\xymatrix@=15mm{
\bbA\ar@/^3mm/@<1mm>@{>->}[r]^{T_{\bbA}}\ar@{}[r]|{\perp} & \P\bbA\ar@/^3mm/@<1mm>@{->>}[r]^F\ar@/^3mm/@<1mm>@{->>}[l]^{\sup_{\bbA}}\ar@{}[r]|{\top} & \R\bbB\ar@/^3mm/@<1mm>@{>->}[l]^G}$$
in $\Cocont(\Q)$, where
$F(\phi)=\Phi\tensor\phi$ and $G(\phi)=\Psi\tensor\phi$.
We can calculate that for $a\in\bbA$, 
\begin{eqnarray*}
\sup_{\bbA}\circ G\circ F\circ T_{\bbA}(a)
 & = & \sup_{\bbA}(\Psi\tensor\Phi\tensor\Theta_{\bbA}(-,a)) \\
 & = & \sup_{\bbA}(\Theta_{\bbA}\tensor\Theta_{\bbA}(-,a)) \\
 & = & \sup_{\bbA}(\Theta_{\bbA}(-,a)) \\
 & \cong & a,
\end{eqnarray*}
using the idempotency of $\Theta_{\bbA}$. For $\phi\in\R\bbB$, it is clear from $T_{\bbA}\dashv\sup_{\bbA}$ that 
$$F\circ T_{\bbA}\circ\sup_{\bbA}\circ G(\phi)\leq F(G(\phi))=\Phi\tensor\Psi\tensor\phi=\bbB\tensor\phi=\phi.$$
For the converse inequality, observe first that
\begin{eqnarray*}
\phi\leq F\circ T_{\bbA}\circ\sup_{\bbA}\circ G(\phi)
 & \iff & G(\phi)\leq T_{\bbA}\circ\sup_{\bbA}\circ G(\phi).
\end{eqnarray*}
But, using that $T_{\bbA}\leq Y_{\bbA}$, we can calculate indeed that
\begin{eqnarray*}
T_{\bbA}\circ\sup_{\bbA}\circ G(\phi)
 & = & \Theta_{\bbA}(-,\sup_{\bbA}\circ G(\phi)) \\
 & = & \Theta_{\bbA}\tensor\bbA(-,\sup_{\bbA}\circ G(\phi)) \\
 & = & \Theta_{\bbA}\tensor\P\bbA(T_{\bbA}-,G(\phi)) \\
 & \geq & \Theta_{\bbA}\tensor\P\bbA(Y_{\bbA}-,G(\phi)) \\
 & = & \Theta_{\bbA}\tensor G(\phi) \\
 & = & G(\phi).
\end{eqnarray*}
This means that $F\circ T_{\bbA}$ and $\sup_{\bbA}\circ G$ constitute the equivalence of $\bbA$ and $\R\bbB$, where $\bbB$ is any regular $\Q$-semicategory over which $\Theta_{\bbA}$ splits. If now $\bbA\simeq\R\bbB'$ for some other regular $\Q$-semicategory $\bbB'$, then $\bbB$ and $\bbB'$ are Morita-equivalent, i.e.\ isomorphic in $\RSDist(\Q)$, so $\Theta_{\bbA}$ also splits over $\bbB'$.
\par
For the converse implication, we've argued above that $\R\bbB$ is totally continuous. And it follows from the first part of the proof that $\Theta_{\R\bbB}$ splits over $\bbB$.
\endofproof
It is an immediate consequence of this important proposition that, for a totally continuous cocomplete $\Q$-category $\bbA$, if $\Theta_{\bbA}\:\bbA\dist\bbA$ splits over some regular $\Q$-semicategory $\bbB$, then $\bbA\simeq\R\bbB$. In particular, 
recalling how idempotents may be split in $\RSDist(\Q)$\footnote{This is a particular case of a general result on the splitting of idempotents in the split-idempotent completion $\Idm(\Q')$ of a given quantaloid $\Q'$, here applied to $\Q'=\Matr(\Q)$. See [Stubbe, 2005b, Appendix] for details.}, we may explicitly say that $\Theta_{\bbA}\:\bbA\dist\bbA$ splits in $\RSDist(\Q)$ over some regular $\Q$-semicategory $\bbB$ if and only if $\bbB$ is Morita equivalent
to the regular $\Q$-semicategory whose $\Q_0$-typed object set is $\bbA_0$, and hom-arrows are $\Theta_{\bbA}(a',a)$ for any $a,a'\in\bbA_0$.

\section{Totally algebraic cocomplete $\Q$-categories}\label{G}

As in section \ref{E}, we write $\Theta_{\bbA}\:\bbA\dist\bbA$ for the totally-below relation on a given cocomplete $\Q$-category $\bbA$ (whether it is totally continuous or not), and the corresponding functor as $T_{\bbA}\:\bbA\to\P\bbA$.
\begin{lemma}\label{400}
Let $\bbA$ be a cocomplete $\Q$-category. For an object $a\in\bbA$, the following are equivalent:
\begin{enumerate}
\item $1_{ta}\leq\Theta_{\bbA}(a,a)$,
\item for all $x\in\bbA$, $\bbA(x,a)\leq\Theta_{\bbA}(x,a)$,
\item for all $x\in\bbA$, $\bbA(a,x)\leq\Theta_{\bbA}(a,x)$,
\item $Y_{\bbA}(a)\leq T_{\bbA}(a)$.
\end{enumerate}
In fact, the ``$\leq$'' may be replaced by ``$=$'' in all statements but the first.
\end{lemma}
\proof
Of course, the second and the fourth statement are tautologies. If $1_{ta}\leq\Theta_{\bbA}(a,a)$ then, for any $x\in\bbA$, $\bbA(x,a)\leq\bbA(x,a)\circ\Theta_{\bbA}(a,a)\leq\Theta_{\bbA}(x,a)$; so the first condition implies the second. Conversely, putting $x=a$ in the second condition, $1_{ta}\leq\bbA(a,a)\leq\Theta_{\bbA}(a,a)$; so the first condition is implied. The equivalence of the first and the third statement is similar. Finally, that the ``$\leq$'' may be replaced by ``$=$'' in statements two to four, is due to \ref{103.0}. 
\endofproof
An object $a\in\bbA$ of a cocomplete $\Q$-category satisfying the equivalent conditions in \ref{400}, is said to be {\em totally compact}. We will write $i\:\bbA\c\to\bbA$ for the full subcategory of $\bbA$ determined by its totally compact objects; it is thus the so-called {\em inverter} of the 2-cell $T_{\bbA}\leq Y_{\bbA}\:\bbA\biar\P\bbA$ in $\Cat(\Q)$, as we spell out next. 
\begin{proposition}\label{401}
For any cocomplete $\Q$-category $\bbA$, the full embedding of the totally compact objects $i\:\bbA\c\to\bbA$ satisfies $T_{\bbA}\circ i\cong Y_{\bbA}\circ i$, and any other functor $F\:\bbC\to\bbA$ such that $T_{\bbA}\circ F\cong Y_{\bbA}\circ F$, factors essentially uniquely through $i$. Moreover, if $F$ is fully faithful, then so is its factorization through $i$.
\end{proposition}
\proof
That $T_{\bbA}\circ i\cong Y_{\bbA}\circ i$, is a rewrite of the fourth condition in \ref{400}. Now assume $T_{\bbA}\circ F\cong Y_{\bbA}\circ F$, i.e.\ $Fc\in\bbA\c$ for each $c\in\bbC$. Since $\bbC(c',c)\leq\bbA(Fc',Fc)=\bbA\c(Fc',Fc)$ already $\overline{F}\:\bbC\to\bbA\c\:c\mapsto Fc$ is a factorization of $F$ through $i$. This factorization is essentially unique, because $i\:\bbA\c\to\bbA$, which is injective on objects, is a monomorphism in $\Cat(\Q)$. It is clear that $\overline{F}$ is fully faithful whenever $F$ is.
\endofproof
It follows straightforwardly that equivalent cocomplete $\Q$-categories, say $\bbA\simeq\bbA'$, have equivalent $\Q$-categories of totally compact objects, $\bbA\c\simeq\bbA\c'$.
\par
For any cocomplete $\Q$-category $\bbA$, we can now {\em define} the distributor $\Sigma_{\bbA}\:\bbA\dist\bbA$ to be precisely the comonad determined by the adjoint pair of distributors induced by the full embedding $i\:\bbA\c\to\bbA$ of totally compact objects:
$$\Sigma_{\bbA}(a',a)=\bbA(a',i-)\tensor\bbA(i-,a).$$
Further we put $S_{\bbA}\:\bbA\to\P\bbA$ to be the functor corresponding to $\Sigma_{\bbA}$ under the universal property of the presheaf category, i.e.\ $S_{\bbA}(a)=\Sigma_{\bbA}(-,a)$.
\begin{lemma}\label{405.0}
For a cocomplete $\Q$-category $\bbA$, $\Sigma_{\bbA}\leq\Theta_{\bbA}$.
\end{lemma}
\proof
By \ref{400} we can calculate that, for any $a,a'\in\bbA$,
$$\Sigma_{\bbA}(a',a)=\bbA(a',i-)\tensor\bbA(i-,a)=\Theta_{\bbA}(a',i-)\tensor\Theta_{\bbA}(i-,a)\leq\Theta_{\bbA}(a',a).$$
This proves our claim.
\endofproof
\par
The following result must be compared with \ref{104}.
\begin{lemma}\label{405}
For a cocomplete $\Q$-category $\bbA$, the following are equivalent:
\begin{enumerate}
\item\label{q} $i\:\bbA\c\to\bbA$ satisfies $\<i,i\>\cong 1_{\bbA}$,
\item\label{r} for every $a\in\bbA$, $\sup_{\bbA}(S_{\bbA}a)\cong a$,
\item\label{s} $\sup_{\bbA}\circ S_{\bbA}\cong 1_{\bbA}$,
\item\label{t} $\colim(\Sigma_{\bbA},1_{\bbA})\cong 1_{\bbA}$.
\end{enumerate}
In this case, $\Sigma_{\bbA}=\Theta_{\bbA}$. 
\end{lemma}
\proof
The left Kan extension of $i\:\bbA\c\to\bbA$ along itself (exists and) is pointwise because $\bbA$ is cocomplete; we may thus compute, using the calculation rules for colimits, that $$\<i,i\>(a)\cong\colim(\bbA(i-,a),i)\cong\sup_{\bbA}(\bbA(-,i-)\tensor\bbA(i,a))=\sup_{\bbA}(S_{\bbA}a).$$ This immediately shows that statements (\ref{q}) and (\ref{r}) are synonymous. It is clear that statements (\ref{r}), (\ref{s}) and (\ref{t}) are synonymous.
\par
Assuming these equivalent conditions to hold, putting $\phi=\Sigma_{\bbA}(-,a)=S_{\bbA}a$ in \ref{104.1} shows that
\begin{eqnarray*}
\Theta_{\bbA}(a',a)
 & \leq & \Big\{\bbA\Big(a,\sup_{\bbA}(S_{\bbA}a)\Big),\Sigma_{\bbA}(a',a)\Big\}\\
 & = & \Big\{\bbA(a,a),\Sigma_{\bbA}(a',a)\Big\} \\
 & \leq & \Sigma_{\bbA}(a',a).
\end{eqnarray*}
But the converse inequality always holds, so we have $\Theta_{\bbA}=\Sigma_{\bbA}$.
\endofproof
Mimicking the classical terminology of [Rosebrugh and Wood, 1994] once more, a cocomplete $\Q$-category is {\em totally algebraic} when it satisfies the equivalent conditions in \ref{405}; that is to say, ``every object is the supremum of the (downclosure of the set of) totally compact objects below it''. 
\par
It is immediate from \ref{405} and \ref{104} that ``totally algebraic'' implies ``totally continuous'', but the converse is not true. (For a counterexample, compare \ref{1000} and \ref{404}, with [Stubbe, 2005b, {\bf 4.7}].)
\begin{proposition}\label{404.0}
For a totally continuous cocomplete $\Q$-category $\bbA$, the following are equivalent:
\begin{enumerate}
\item $\bbA$ is totally algebraic,
\item $\bbA$ is totally continuous and $\Theta_{\bbA}=\Sigma_{\bbA}$.
\end{enumerate}
\end{proposition}
\proof
For the non-trivial implication, note that $$\colim(\Sigma_{\bbA},1_{\bbA})=\colim(\Theta_{\bbA},1_{\bbA})\cong 1_{\bbA}$$ 
whenever $\bbA$ is totally continuous and $\Theta_{\bbA}=\Sigma_{\bbA}$.
\endofproof
The following should be compared with \ref{1000}.
\begin{proposition}\label{404}
For a totally continuous cocomplete $\Q$-category $\bbA$, the following are equivalent:
\begin{enumerate}
\item $\bbA$ is totally algebraic,
\item $\bbA\simeq\P\bbA\c$,
\item $\bbA\simeq\P\bbC$ for some $\Q$-category $\bbC$.
\end{enumerate}
\end{proposition}
\proof
It follows directly from \ref{405} that for a totally algebraic $\bbA$, $\Theta_{\bbA}(=\Sigma_{\bbA})$ splits over the $\Q$-category $\bbA\c$; so \ref{1000} implies that $\bbA\simeq\P\bbA\c$.
\par
Suppose now that $\bbA\simeq\P\bbC$ for some $\Q$-category $\bbC$; by \ref{1000} we know that $\bbA$ is totally continuous and that there is a splitting
$$\xymatrix@=15mm{
\bbA\ar@(u,l)_{\Theta_{\bbA}}|{\distsign}\ar@/^3mm/[r]^{\Phi}|{\distsign} & \bbC\ar@/^3mm/[l]^{\Psi}|{\distsign}}$$
of the comonad $\Theta_{\bbA}$ in $\Dist(\Q)$. Then in particular $\Psi\dashv\Phi$, and therefore -- since any cocomplete $\Q$-category is Cauchy complete\footnote{See [Stubbe, 2005a, section 7] for a presentation of the theory of Cauchy complete $\Q$-categories.} -- there exists a functor $F\:\bbC\to\bbA$ such that $\Psi=\bbA(-,F-)$ and $\Phi=\bbA(F-,-)$. Observe that for each $c\in\bbC$, 
$$\Theta_{\bbA}(Fc,Fc)=\bbA(Fc,F-)\tensor\bbA(F-,Fc)\geq 1_{tFc},$$
i.e.\ each $Fc$ is totally compact in $\bbA$. So $F$ factors over $i\:\bbA\c\to\bbA$ by some functor $\overline{F}\:\bbC\to\bbA\c$ (cf.~\ref{401}), and we have that
\begin{eqnarray*}
\Theta_{\bbA}(-,-)
 & = & \bbA(-,F-)\tensor\bbA(F-,-) \\
 & = & \bbA(-,i-)\tensor\bbA\c(-,\overline{F}-)\tensor\bbA\c(\overline{F}-,-)\tensor\bbA(i-,-) \\
 & \leq & \bbA(-,i-)\tensor\bbA(i-,-) \\
 & = & \Sigma_{\bbA}.
\end{eqnarray*}
So we conclude that $\Theta_{\bbA}=\Sigma_{\bbA}$ (because the converse inequality always holds) and, by \ref{404.0}, $\bbA$ is totally algebraic.
\endofproof
From this proof it follows that a cocomplete $\bbA$ is totally algebraic if and only if there exist a $\Q$-category $\bbC$ and a fully faithful functor $F\:\bbC\to\bbA$ such that $\Theta_{\bbA}$ is the comonad determined by the adjunction $\bbA(-,F-)\dashv\bbA(F-,-)$ in $\Dist(\Q)$; and that in this case {\em every} splitting of $\Theta_{\bbA}$ in $\Dist(\Q)$ is of this kind.

\section{Cauchy completions revisited}\label{H}

Already in the proof of \ref{404}, the theory of Cauchy complete $\Q$-categories comes lurking around the corner. We can exhibit a more explicit link.
\par
First observe that from \ref{1000} we know that, for any $\Q$-category $\bbC$, the presheaf category $\P\bbC$ is totally continuous and that the totally-below relation $\Theta_{\P\bbC}$ splits over $\bbC$; and from \ref{404} we know that $\P\bbC$ is even totally algebraic and that there must be a fully faithful functor $F\:\bbC\to\P\bbC$ such that $\Theta_{\P\bbC}$ is the comonad determined by the adjunction $\P\bbC(-,F-)\dashv\P\bbC(F-,-)$. The following lemma shows that it is the Yoneda embedding $Y_{\bbC}\:\bbC\to\P\bbC$ that does the job.
\begin{lemma}\label{402}
For a $\Q$-category $\bbC$, the totally-below relation on $\P\bbC$ is $$\Theta_{\P\bbC}=\P\bbC(-,Y_{\bbC}-)\tensor\P\bbC(Y_{\bbC}-,-).$$
\end{lemma}
\proof
The fully faithful Yoneda embedding $Y_{\bbC}\:\bbC\to\P\bbC$ induces an adjoint pair 
$$\xymatrix@=20mm{
\bbC\ar@{}[r]|{\perp}\ar@/^3mm/[r]|{\distsign}^{\P\bbC(-,Y_{\bbC}-)}
 & \P\bbC\ar@/^3mm/[l]|{\distsign}^{\P\bbC(Y_{\bbC}-,-)}}$$
in $\Dist(\Q)$, the unit of the adjunction being an equality. Applying the local equivalence (\ref{2}) gives
$$\xymatrix@=15mm{
\P\bbC\ar@{}[r]|{\perp}\ar@{>->}@/^3mm/[r]^{F} & \P\P\bbC\ar@{->>}@/^3mm/[l]^{G}}$$
in $\Cocont(\Q)$, where $F(\phi)=\P\bbC(-,Y_{\bbC}-)\tensor\phi$ and $G(\Phi)=\P\bbC(Y_{\bbC}-,-)\tensor\Phi$ for $\phi\in\P\bbC$ and $\Phi\in\P\P\bbC$. But the calculation rules for colimits in presheaf categories [Stubbe, 2005a, {\bf 6.4}] imply that $G=\sup_{\P\bbC}$. This means that its left adjoint $F$ is actually $T_{\P\bbC}$, and thus that, for $\phi,\psi\in\P\bbC$, $$\Theta_{\P\bbC}(\psi,\phi)=F(\phi)(\psi)=\P\bbC(\psi,Y_{\bbC}-)\tensor\phi=\P\bbC(\psi,Y_{\bbC}-)\tensor\P\bbC(Y_{\bbC}-,\phi).$$
This proves our claim.
\endofproof
\begin{proposition}\label{402.1}
For a $\Q$-category $\bbC$, the category $(\P\bbC)\c$ of totally compact objects in $\P\bbC$ is (equivalent to) the Cauchy completion $\bbC\cc$ of $\bbC$.
\end{proposition}
\proof
We will show that a presheaf $\phi\in\P\bbC$ is a Cauchy presheaf (i.e.\ that it has a right adjoint in $\Dist(\Q)$) if and only if it is a totally compact object in $\P\bbC$. First assume that $\phi$ is totally compact; using \ref{402} this means that $1_{t\phi}\leq\P\bbC(\phi,Y_{\bbC}-)\tensor\phi$. But is also true that $\phi\tensor\P\bbC(\phi,Y_{\bbC}-)=\P\bbC(Y_{\bbC}-,\phi)\tensor\P\bbC(\phi,Y_{\bbC}-)\leq\P\bbC(Y_{\bbC}-,Y_{\bbC}-)=\bbC$. This proves that $\phi\dashv\P\bbC(\phi,Y_{\bbC}-)$. Conversely, suppose that $\phi\dashv\phi^*$, then necessarily\footnote{If, in a locally ordered category $\K$, a morphism $f\:A\to B$ is known to have a right adjoint $f^*$, then the right lifting $[f,1_B]$ exists and equals $f^*$.} $\phi^*=[\phi,\bbC]=\P\bbC(\phi,Y_{\bbC}-)$, and therefore $1_{t\phi}\leq\phi^*\tensor\phi=\P\bbC(\phi,Y_{\bbC}-)\tensor\P\bbC(Y_{\bbC}-,\phi)$. By \ref{402} this means that $\phi$ is totally compact in $\P\bbC$. To conclude the proof, recall that the full subcategory of $\P\bbC$ determined by the Cauchy presheaves is indeed the Cauchy completion $\bbC\cc$ of $\bbC$.
\endofproof
It follows now from \ref{404} and \ref{402.1} that for a totally algebraic cocomplete $\Q$-category $\bbA$, the full subcategory $\bbA\c$ of totally compact objects is Cauchy complete: because $\bbA\simeq\P\bbC$ implies $\bbA\c\simeq(\P\bbC)\c\simeq\bbC\cc$, and a category which is equivalent to a Cauchy complete category is Cauchy complete itself.

\section{In terms of modules}\label{I}

The locally ordered category $\Cocont(\Q)$ is biequivalent to $\QUANT(\Q\op,\Sup)$, the quantaloid of (right) $\Q$-modules. Explicitly, a cocomplete $\Q$-category $\bbA$ determines the module
$${\cal F}_{\bbC}\:\Q\op\to\Sup\:\Big(f\:X\to Y\Big)\mapsto\Big(-\tensor f\:\bbC_Y\to\bbC_X\Big).$$
And a module ${\cal F}\:\Q\op\to\Sup$ determines the cocomplete $\Q$-category $\bbC_{\cal F}$ with objects $\coprod_X{\cal F}X$ (and $tx=X$ if and only if $x\in{\cal F}X$) and hom-arrows 
$$\bbC_{\cal F}(y,x)=\bigvee\{f\:ty\to tx\mid {\cal F}(f)(y)\leq x\}.$$ 
This is really a part of the theory of tensored and cotensored $\Q$-categories; [Stubbe, 2004, section 4] contains the details. It is then a matter of fact that the projective objects in $\Cocont(\Q)$ correspond to those in $\QUANT(\Q\op,\Sup)$ under this biequivalence. 
\begin{proposition}\label{a0}
Let $\bbA$ and $\cal F$ be a cocomplete $\Q$-category and a $\Q$-module that correspond to each other under the biequivalence 
$$\Cocont(\Q)\simeq\QUANT(\Q\op,\Sup)$$ then the following are equivalent:
\begin{enumerate}
\item $\bbA$ is a projective object of $\Cocont(\Q)$,
\item $\cal F$ is a projective object of $\QUANT(\Q\op,\Sup)$.
\end{enumerate}
\end{proposition}
\par
Since $\QUANT(\Q\op,\Sup)$ is a (large) quantaloid (in particular -- and in contast to $\Cocont(\Q)$ -- its local order is anti-symmetric), an object ${\cal F}$ is projective if and only if the representable homomorphism 
\begin{equation}\label{a6}
\QUANT(\Q\op,\Sup)({\cal F},-)\:\QUANT(\Q\op,\Sup)\to\Sup
\end{equation}
preserves epimorphisms. (This is really a straightforward reformulation of the definition of ``projectivity'' that was given in section \ref{C}.) A seemingly stronger notion is of much importance in the theory of ($\Sup$-)enriched categories: after [Kelly, 1982], a {\em small-projective object} ${\cal F}\in\QUANT(\Q\op,\Sup)$ is one for which the representable homomorphism in (\ref{a6}) preserves all small weighted colimits. Clearly a small-projective object in $\QUANT(\Q\op,\Sup)$ is also projective; but we will prove that the converse also holds. Thereto we exploit the biequivalence between $\Q$-modules and cocomplete $\Q$-categories, using the notion of a ``truly free object''. Part of this stems from [van der Plancke, 1997].
\par
There is a forgetful functor\footnote{This is not a 2-functor, for it is not defined on 2-cells; so the adjunction it is part of, is not a 2-adjunction!} $(-)_0\:\Cat(\Q)\to\Set/\Q_0$ sending a $\Q$-category $\bbA$ to its underlying $\Q_0$-typed set of objects. This forgetful admits a left adjoint: it sends a $\Q_0$-typed set $A$ to the ``identity matrix'' on $A$; we denote it $A\f$. The unit of this adjunction is the identity. For an $\bbA\in\Cat(\Q)$, the component at $\bbA$ of the counit of this adjunction is the functor $\bbA_0\f\to\bbA\:a\mapsto a$ (which is the identity on objects, but not on hom-arrows!). By composition of the adjunctions
$$\Cocont(\Q)
\xymatrix@=15mm{ \ar@{}[r]|{\perp}\ar@/_3mm/[r]_{\U} & \ar@/_3mm/[l]_{\P} } 
\Cat(\Q) 
\xymatrix@=15mm{ \ar@{}[r]|{\perp}\ar@/_3mm/[r]_{(-)_0} & \ar@/_3mm/[l]_{(-)\f} } 
\Set/\Q_0$$
we obtain a ``truly forgetful'' functor---and by a {\em truly free object} in $\Cocont(\Q)$ we will mean a free object relative to this truly forgetful functor, i.e.\ an object equivalent to $\P A\f$ for some $\Q_0$-typed set $A$. The component at $\bbA\in\Cocont(\Q)$ of the counit of $\P\circ(-)\f\dashv(-)_0\circ\U$ is the cocontinuous functor 
$$\xymatrix@=15mm{
\P\bbA_0\f\ar[r]^{\P(\varepsilon^1_{\U(\bbA)})} & \P\bbA\ar[r]^{\varepsilon^2_{\bbA}} & \bbA}$$
where $\varepsilon^1$ and $\varepsilon^2$ are the counits of, respectively, $(-)\f\dashv(-)_0$ and $\P\dashv\U$. We already know from section \ref{C} that $\varepsilon^2_{\bbA}=\sup_{\bbA}$ is surjective. A straightforward calculation shows that also $\P(\varepsilon^1_{\U(\bbA)})$ is surjective. The following must now be compared with \ref{7}.
\begin{proposition}\label{a1}
For a cocomplete $\Q$-category $\bbA$, the following are equivalent:
\begin{enumerate}
\item $\bbA$ is a projective object in $\Cocont(\Q)$,
\item $\bbA$ is a retract of $\P\bbA_0\f$ in $\Cocont(\Q)$,
\item $\bbA$ is a retract of a truly free object in $\Cocont(\Q)$.
\end{enumerate}
\end{proposition}
\proof
If $\bbA$ is projective, then $1_{\bbA}$ factors through the composite surjection $\P\bbA_0\f\surj\P\bbA\surj\bbA$ -- the component at $\bbA$ of the counit of the adjunction explained above -- so that $\bbA$ is a retract of $\P\bbA_0\f$. Obviously, the second statement implies the third. And a truly free object is free, so the third statement implies the first (by \ref{7} for example). 
\endofproof
\par
Now we will translate the equivalence of the first and the third statement in \ref{a1} from $\Cocont(\Q)$ to the biequivalent $\QUANT(\Q\op,\Sup)$.
\begin{lemma}\label{a21}
Every truly free object of $\Cocont(\Q)$ is the coproduct of truly free objects on singletons.
\end{lemma}
\proof
Since left adjoints preserve coproducts and each $\Q_0$-typed set $A$ is (in the obvious way) the coproduct of singletons, it follows that $\P A\f=\coprod_{a\in A}\P(\{a\}\f)$. 
\endofproof
A singleton object of $\Set/\Q_0$ is, essentially, a ``duplicate'' of an object of $\Q$: a singleton $\{a\}\in\Set/\Q_0$ determines the object $ta\in\Q$, and an object $X\in\Q$ determines the singleton $\{*\}\in\Set/\Q_0$ whose single object is of type $X$. This correspondence is essentially bijective. In [Stubbe, 2004, 2005a, 2005b] 
we have, for a given object $X\in\Q$, systematically denoted $*_X$ for the free $\Q$-category on the singleton determined by $X$; and $\P X$ was our notation for the presheaf category on such a $*_X$. That is to say, those $\{\P X\mid X\in \Q\}$ are essentially the truly free objects on singletons of $\Cocont(\Q)$.
\begin{lemma}\label{a2}
The truly free objects on singletons of $\Cocont(\Q)$ correspond under the biequivalence with $\QUANT(\Q\op,\Sup)$ to representable modules.
\end{lemma}
\proof
Given an ${\cal F}=\Q(-,X)\:\Q\op\to\Sup$ it is easily verified that the $\Q$-category $\bbC_{\cal F}$ is $\P X$. Conversely, for a $\P X$, with $X\in\Q$, it is easily seen that the module ${\cal F}_{\bbC}$ is represented by $X\in\Q$. This correspondence is bijective.
\endofproof
Here is, then, the conclusion to the previous lemmas.
\begin{proposition}\label{a3}
The projective objects of $\QUANT(\Q\op,\Sup)$ are precisely the retracts of direct sums of representable modules.
\end{proposition}
\par
Finally we make the link with small-projectives in $\QUANT(\Q\op,\Sup)$. It is proved in [Kelly, 1982, {\bf 5.26}] (in the more general context of $\V$-enriched categories) that representable $\Q$-modules are small-projective; and [Kelly,  1982, {\bf 5.25}] shows that retracts of small-projective $\Q$-modules are small-projective themselves. In the specific case of $\Sup$-enrichment, using that in any quantaloid sums and products coincide, we may also prove the following.
\begin{lemma}\label{a7}
A direct sum of small-projective $\Q$-modules is small-projective.
\end{lemma}
\proof
Consider a (set-indexed) family $(\F_i)_{i\in I}$ of small-projective $\Q$-modules, and a small weighted colimit diagram
$$\xymatrix@=15mm{
{\cal I}\ar[r]|{\distsign}^{\Phi} & {\cal R}\ar[r]^D & }\QUANT(\Q\op,\Sup).$$
As is customary, ${\cal I}$ stands for the one-object quantaloid whose hom-object is the identity for the tensor in $\Sup$. We may then calculate in $\Sup$ that
\begin{eqnarray*}
 & & \QUANT(\Q\op,\Sup)\Big(\oplus_{i\in I}\F_i,\colim(\Phi,D)\Big) \\
 & & \cong \oplus_{i\in I}\enskip\QUANT(\Q\op,\Sup)\Big(\F_i,\colim(\Phi,D)\Big) \\
 & & \cong \oplus_{i\in I}\enskip\colim\Big(\Phi,\QUANT(\Q\op,\Sup)(\F_i,D-)\Big) \\
 & & \cong \colim\Big(\Phi,\oplus_{i\in I}\QUANT(\Q\op,\Sup)(\F_i,D-)\Big) \\
 & & \cong \colim\Big(\Phi,\QUANT(\Q\op,\Sup)\Big(\oplus_{i\in I}\F_i,D-\Big)\Big).
\end{eqnarray*}
The contravariant representable represented by $\colim(\Phi,D)$ turns sums into products (but both are direct sums). Then we use the hypothetical small-projectivity of the $\F_i$ and the ``general interchange of colimits'' [Kelly, 1982, {\bf (3.21)}]. Finally the contravariant representable represented by $D-$ turns products into sums (but both are direct sums).
\endofproof
Because a small-projective is always projective, \ref{a3}, \ref{a7} and the theorems in [Kelly, 1982] recalled above, imply the following.
\begin{proposition}\label{a8}
For ${\cal F}\in\QUANT(\Q\op,\Sup)$, the following are equivalent:
\begin{enumerate}
\item ${\cal F}$ is a projective object,
\item ${\cal F}$ is a retract of a direct sum of representable $\Q$-modules,
\item ${\cal F}$ is a small-projective object.
\end{enumerate}
\end{proposition}
Via \ref{a0} this says something about projective objects in $\Cocont(\Q)$ too.



\begin{thebibliography}{10}

\bibitem{abramskyandjung}
Abramsky, Samson and Achim Jung, {\em Domain theory}, Handbook of logic in computer science (volume 3), Oxford University Press (1994), pp.~1--168.

\bibitem{abramskyandvickers}
Abramsky, Samson and Steve Vickers, {\em Quantales, observational logic and process semantics}, Math. Struct. in Comput. Sci. \textbf{3} (1993), pp.~161--227.

\bibitem{adamek}
Ad{\'a}mek, Ji{\v{r}}{\'{\i}}, {\em A categorical generalization of Scott domains}, Math. Struct. in Comput. Sci. \textbf{7} (1997), pp.~419-443.

\bibitem{adamekandrosicky}
Ad{\'a}mek, Ji{\v{r}}{\'{\i}} and Ji{\v{r}}{\'{\i}} Rosick{\'y}, ``Locally presentable and accessible categories'', Cambridge University Press, Cambridge, 1994.

\bibitem{baltagandcoeckeandsadrzadeh} Baltag, Alexandru, Bob Coecke and Mehrnoosh Sadrzadeh, {\em Epistemic actions as resources}, paper presented at the workshop on Logics for Resources, Processes, and Programs, Turku, Finland, 2004.

\bibitem{borceux}
Borceux, Francis, ``Handbook of categorical algebra (3 volumes)'', Cambridge University Press, Cambridge, 1994.

\bibitem{coeckemooreandstubbe}
Coecke, Bob, David Moore and Isar Stubbe, {\em Quantaloids describing causation and propagation of physical properties}, Found. Phys. Lett. \textbf{14} (2001), pp.~133--145. 

\bibitem{coeckeandstubbe}
Coecke, Bob and Isar Stubbe, {\em Operational resolutions and state transitions in a categorical setting}, Found. Phys. Letters \textbf{12} (1999), pp.~29--49. 

\bibitem{fawcettandwood}
Fawcett, Barry and Richard~J. Wood, {\em Constructive complete distributivity I}, Math. Proc. Cam. Phil. Soc. \textbf{107} (1990), pp.~81--89.

\bibitem{floyd}
Floyd, Robert, {\em Assigning meanings to programs}, Proc. Symp. Applied Math. \textbf{19} (1967), pp.~19--32.

\bibitem{gordonandpower}
Gordon, Robert and A. John Power, {\em Enrichment through variation}, J. Pure Appl. Algebra \textbf{120} (1997), pp.~167--185.

\bibitem{hoare}
Hoare, Tony, {\em An Axiomatic Basis for Computer Programming}, Comm. of the ACM \textbf{12} (1969), pp.~576--583.

\bibitem{kelly}
Kelly, G. Max, ``Basic concepts of enriched category theory'', Cambridge University Press, Cambridge, 1982. Also in: Reprints in Theory Appl. of Categ. \textbf{10} (2005), pp.~1--136.
  
\bibitem{kellyandschmitt}
Kelly, G. Max and Vincent Schmitt, {\em Notes on enriched categories with colimits of some class}, preprint (2005), {\tt arXiv:math.CT/0509102}.

\bibitem{lawvere} 
Lawvere, Bill, {\em Metric spaces, generalized logic and closed categories}, Rend. Sem. Mat. Fis. Milano \textbf{43} (1973), pp.~135--166. Also in: Reprints in Theory Appl. of Categ. \textbf{1} (2002), pp.~1--37.

\bibitem{resende} 
Resende, Pedro, {\em Quantales and observational semantics}, Fund. Theories Phys. \textbf{111} (``Current research in operational quantum logic''), Kluwer Academic Publishers (2000), pp.~263--288.

\bibitem{rosebrughandwood:06030}
Rosebrugh, Robert and Richard~J. Wood, {\em Constructive complete distributivity IV}, Appl. Categ. Structures \textbf{2} (1994), pp.~119--144.

\bibitem{isar} 
Stubbe, Isar, {\em Causal duality: what it is, and what it is good for}, paper presented at the Sixth Biannual Meeting of the International Quantum Structures Association, Vienna, Austria, 2002.

\bibitem{isarc} 
Stubbe, Isar, {\em Categorical structures
enriched in a quantaloid: tensored and cotensored categories}, preprint (2004), {\tt arXiv:math.CT/0411366}.
  
\bibitem{isara} 
Stubbe, Isar, {\em Categorical structures
enriched in a quantaloid: categories, distributors and functors}, Theory Appl. Categ. \textbf{14} (2005a), pp.~1--45.

\bibitem{isarb} 
Stubbe, Isar, {\em Categorical structures
enriched in a quantaloid: regular presheaves, regular semicategories}, Cahiers Topol. G\'eom. Diff\'er. Cat\'eg. \textbf{46} (2005b), pp.~99--121.

\bibitem{isard}
Stubbe, Isar, {\em Towards ``dynamical domains'': totally continuous cocomplete $\Q$-categories (Extended abstract)}, Elec. Notes Theor. Comp. Sc. (2006), to appear.

\bibitem{vanderplancke}
van~der Plancke, Fr\'ed\'eric, ``Sheaves on a quantaloid as
  enriched categories without units'', PhD thesis, Universit\'e de Louvain, Louvain-la-Neuve, 1997.

\bibitem{wagner} 
Wagner, Kim R., {\em Liminf convergence in $\Omega$-categories}, Theor. Comput. Science \textbf{184} (1997), pp.~61--104.

\bibitem{wood}
Wood, Richard J., {\em Ordered sets via adjunctions}, Categorical foundations, Cambridge University Press (2004), pp.~5--47.

\end{thebibliography}
\end{document}